\documentclass[12pt]{article}%
\usepackage{amsfonts}
\usepackage{sw20bams}
\usepackage{amsmath}
\usepackage{amssymb}
\usepackage{graphicx}%
\providecommand{\U}[1]{\protect\rule{.1in}{.1in}}
\begin{document}

\title{Half-Iterates and Delta Conjectures}
\author{Steven Finch}
\date{September 3, 2025}
\maketitle

\begin{abstract}
The vivid contrast between two competing algorithms for solving Abel's
equation $g(\theta(x))=g(x)+1$, given $\theta(x)$, is easily sketched.
\ EJ\ is faster and more efficient, but ML evaluates a limit characterizing
the principal solution $g(x)$ directly. \ EJ\ finds $g(x)+\delta$, where
$\delta$ is possibly nonzero but independent of $x$. \ If we were to know an
exact expression for $\delta$, then the \textquotedblleft
intrinsicality\textquotedblright\ of ML would be subsumed by EJ. \ Filling
this gap in our knowledge is the aim of this paper.

\end{abstract}

\footnotetext{Copyright \copyright \ 2025 by Steven R. Finch. All rights
reserved.}

The Mavecha-Laohakosol \cite{dB-zeal, BR-zeal, ML-zeal, F0-zeal} and
\'{E}calle-Jagy \cite{E1-zeal, Kc-zeal, F1-zeal, F2-zeal} algorithms enable us
to numerically solve functional equations%
\[%
\begin{array}
[c]{ccc}%
g_{13}\left(  x\left(  1-x^{2}\right)  \right)  =g_{13}(x)+1, &  & 0<x<1;
\end{array}
\]%
\[%
\begin{array}
[c]{ccc}%
g_{14}\left(  x(1-x+x^{2})\right)  =g_{14}(x)+1, &  & 0<x<1;
\end{array}
\]%
\[%
\begin{array}
[c]{ccc}%
g_{15}\left(  x(1-x)^{2}\right)  =g_{15}(x)+1, &  & 0<x<2;
\end{array}
\]%
\[%
\begin{array}
[c]{ccc}%
g_{16}\left(  x(1-3x+3x^{2})\right)  =g_{16}(x)+1, &  & 0<x<1
\end{array}
\]
in highly distinct ways. \ It is useful to have both methods available. \ For
both, the base function $\theta(x)$ is assumed to be analytic with values
$\theta(0)=0$ \& $\theta^{\prime}(0)=1$; further, its Taylor series at the
origin is%
\[%
\begin{array}
[c]{ccccc}%
x+%
{\displaystyle\sum\limits_{m=1}^{\infty}}
c_{m}x^{m\,\tau+1}, &  & c_{1}=\gamma<0, &  & \tau\geq1\text{ is an integer.}%
\end{array}
\]
In words, $\theta(x)$ possesses an attractive fixed point at $x=0$ with unit
slope and, beyond this, series terms at uniform separation $\tau$. \ Theory
suggests that these conditions are sufficient for the ML\ algorithm to
converge. \ It seems that the same is true for the EJ\ algorithm, which is
simpler than ML and outperforms ML\ in every regard.

Our approach is experimental, focusing mostly on cubics. \ ML computes the
principal solution $g(x)$ of Abel's equation directly, whereas the
EJ\ algorithm constructs $\lambda(x)$ that solves Julia's equation%
\[
\lambda(\theta(x))=\theta^{\prime}(x)\lambda(x)
\]
whose reciprocal approximates the derivative $g^{\prime}(x)$. \ EJ\ has a
drawback: it provides not $g(x)$ but instead $g(x)+\delta$, where $\delta$ is
(usually) a nonzero constant. \ Let\ $\theta_{13}(x)=x\left(  1-x^{2}\right)
$. \ We wish to conjecture an exact expression for $\delta_{13}$, among other
things. \ Clearly $\tau=2$ and $\gamma=-1$. \ EJ\ gives
\begin{align*}
g_{13}(x)  &  =\frac{1}{2x^{2}}+\frac{3}{2}\ln(x)+\frac{5}{8}x^{2}+\frac
{21}{32}x^{4}+\frac{35}{32}x^{6}+\frac{2717}{1280}x^{8}+\frac{13429}%
{3200}x^{10}+\frac{81239}{10752}x^{12}\\
&  +\frac{1271651}{125440}x^{14}+\frac{521379}{286720}x^{16}-\frac
{18534919}{460800}x^{18}-\frac{1204285773}{11264000}x^{20}\\
&  +\frac{1901417927}{20815872}x^{22}+\frac{2501877529163}{1640038400}%
x^{24}+\cdots.
\end{align*}
To actually calculate $g_{13}(x)$ for some $0<x<1$, define $x_{0}=x$ and
$x_{n}=x_{n-1}\left(  1-x_{n-1}^{2}\right)  $ for $n\geq1$. \ From the
definition of $g_{13}(x)$, we have
\[
g_{13}(x)=\lim_{n\rightarrow\infty}\,\left[  g_{13}(x_{n})-n\right]  .
\]
For example, to 100 decimal digits of accuracy,%
\begin{align*}
g_{13}\left(  \frac{1}{2}\right)   &
=1.18935108880763516558234924956664105059667671443624\backslash\\
&  \;\;\;\;\;\;\;\;71933977369986198458248324365808279005576819103340...,
\end{align*}%
\begin{align*}
g_{13}\left(  \frac{1}{\sqrt{3}}\right)   &
=1.05481186044739525804034624653292094298561372154532\backslash\\
&  \;\;\;\;\;\;\;\;98153825514557676637413068364378707969808390093840....
\end{align*}
ML\ allows calculation of%
\[
\tilde{g}_{13}(x)=-2\lim_{n\rightarrow\infty}n^{3/2}\left(  \sqrt{2}%
\,x_{n}-\frac{1}{n^{1/2}}+\frac{3}{8}\frac{\ln(n)}{n^{3/2}}\right)
\]
to high precision and \cite{F3-zeal}
\begin{align*}
\tilde{g}_{13}\left(  \frac{1}{2}\right)   &
=1.70921147422759414764527334066027347665330181520643\backslash\\
&  \;\;\;\;\;\;\;\;86339882470057398910413097076175322980529292243496...,
\end{align*}%
\begin{align*}
\tilde{g}_{13}\left(  \frac{1}{\sqrt{3}}\right)   &
=1.57467224586735424010327033762655336904223882231552\backslash\\
&  \;\;\;\;\;\;\;\;12559730614628877089577841074745751944760863233997....
\end{align*}
The difference between the EJ-based \&\ ML-based values is nonzero and
conjectured to be%
\[
\delta_{13}=g_{13}\left(  1/2\right)  -\tilde{g}_{13}\left(  1/2\right)
=g_{13}\left(  1/\sqrt{3}\right)  -\tilde{g}_{13}\left(  1/\sqrt{3}\right)
=-(3/4)\ln(2).
\]
This identity should hold for any choice of $x=x_{0}$. \ A graph of $\tilde
{g}_{13}\left(  x\right)  $ appears in Figure 1. \ The half-iterate of
$\theta_{13}(x)$ is defined by%
\[%
\begin{array}
[c]{ccc}%
\theta_{13}^{[1/2]}(x)=g_{13}^{[-1]}\left(  g_{13}(x)+\dfrac{1}{2}\right)  &
& \text{if }0<x<1
\end{array}
\]
and is plotted in Figure 2. \ For example,
\begin{align*}
\theta_{13}^{[1/2]}\left(  \frac{1}{2}\right)   &
=0.41820429029068859160564663263510996432355959664258\backslash\\
&  \;\;\;\;\;\;\;\;55277260336535008451118441673117069901712981511710...,
\end{align*}%
\begin{align*}
\theta_{13}^{[1/2]}\left(  \frac{1}{\sqrt{3}}\right)   &
=0.43400787521828089000086877815934526729705657496748\backslash\\
&  \;\;\;\;\;\;\;\;10398836584355586827352894865816690477408133034496....
\end{align*}
Generalizations include%
\[
\delta(k)=-\frac{k+1}{2k}\ln(k)
\]
corresponding to $\theta(x)=x\left(  1-x^{k}\right)  $ for integer $k\geq1$
and
\[
\delta(\rho)=\frac{1-\rho+\rho^{2}}{(1-\rho)^{2}}\ln\left(  \frac{1}{1-\rho
}\right)
\]
corresponding to $\theta(x)=x\left(  1-x\right)  (1+\rho\,x)$ for $0<\rho<1$.
\ For the former,%
\[
\tilde{g}(x)=-k\lim_{n\rightarrow\infty}n^{(k+1)/k}\left(  k^{1/k}%
\,x_{n}-\frac{1}{n^{1/k}}+\frac{k+1}{2k^{2}}\frac{\ln(n)}{n^{(k+1)/k}}\right)
\]
which resembles the case $k=2$ listed earlier. \ For the latter,%
\[
\tilde{g}(x)=-(1-\rho)\lim_{n\rightarrow\infty}n^{2}\left(  \,x_{n}-\frac
{1}{(1-\rho)n}+\frac{1-\rho+\rho^{2}}{(1-\rho)^{3}}\frac{\ln(n)}{n^{2}%
}\right)
\]
which does \textit{not} resemble the case $\rho=1$ listed earlier. \ Various
other delta conjectures are listed in \cite{F2-zeal}, building on work in
\cite{F1-zeal}. \ 

We now turn to $\theta_{14}(x)=x(1-x+x^{2})$, for which $\tau=1$ and
$\gamma=-1$. \ EJ\ gives%
\begin{align*}
g_{14}(x)  &  =\frac{1}{x}-x-x^{2}-x^{3}-\frac{3}{4}x^{4}-\frac{1}{10}%
x^{5}+\frac{3}{4}x^{6}+\frac{13}{14}x^{7}-\frac{3}{4}x^{8}-\frac{167}{45}%
x^{9}\\
&  -\frac{57}{20}x^{10}+\frac{5081}{660}x^{11}+\frac{1837}{120}x^{12}%
-\frac{24235}{1092}x^{13}-\frac{14489}{168}x^{14}+\frac{269263}{3150}%
x^{15}+\cdots.
\end{align*}
To assess $g_{14}(x)$ for some $0<x<1$, define $x_{0}=x$ and $x_{n}%
=x_{n-1}\left(  1-x_{n-1}+x_{n-1}^{2}\right)  $ for $n\geq1$. \ ML\ allows
calculation of%
\[
\tilde{g}_{14}(x)=-\lim_{n\rightarrow\infty}n^{2}\left(  \,x_{n}-\frac{1}%
{n}\right)
\]
to high precision and
\begin{align*}
\tilde{g}_{14}\left(  \frac{1}{2}\right)   &
=1.08497872758210320613109766431163045230197877946729\backslash\\
&  \;\;\;\;\;\;\;\;58035039848467035694746977730608826554708937283892...,
\end{align*}%
\begin{align*}
\tilde{g}_{14}\left(  \frac{1}{3}\right)   &
=2.51000182711894064597232155164267425915580827922384\backslash\\
&  \;\;\;\;\;\;\;\;06758952435313433766715179785405812086576884568208....
\end{align*}
We conjecture that $\delta_{14}=0$ and this holds for any choice of $x=x_{0}$.
\ A graph similar to that of $\tilde{g}_{14}\left(  x\right)  $ appears in
\cite{F4-zeal}, associated with $x(1-a\,x+x^{2})$ for $a=3/2$ (as opposed to
$a=1$). \ A\ generalization is $\delta(k)=0$ corresponding to both%
\[%
\begin{array}
[c]{ccccc}%
x\,%
{\displaystyle\sum\limits_{j=0}^{k}}
(-1)^{j}x^{j} &  & \text{and} &  & \dfrac{x}{1+x}%
\end{array}
\]
for integer $k\geq1$; the limiting case $x/(1+x)$ is trivial \cite{F3-zeal}.

We now turn to $\theta_{15}(x)=x(1-x)^{2}$, for which $\tau=1$ and $\gamma
=-2$. \ EJ\ gives
\begin{align*}
g_{15}(x)  &  =\frac{1}{2x}+\frac{3}{4}\ln(x)+\frac{5}{8}x+\frac{21}{32}%
x^{2}+\frac{35}{32}x^{3}+\frac{2717}{1280}x^{4}+\frac{13429}{3200}x^{5}%
+\frac{81239}{10752}x^{6}+\frac{1271651}{125440}x^{7}\\
&  +\frac{521379}{286720}x^{8}-\frac{18534919}{460800}x^{9}-\frac
{1204285773}{11264000}x^{10}+\frac{1901417927}{20815872}x^{11}+\frac
{2501877529163}{1640038400}x^{12}\\
&  +\frac{779013297328117}{223865241600}x^{13}-\frac{10607531187899069}%
{803618816000}x^{14}-\frac{48418902026319109}{469647360000}x^{15}-\cdots.
\end{align*}
To assess $g_{15}(x)$ for some $0<x<2$, define $x_{0}=x$ and $x_{n}%
=x_{n-1}\left(  1-x_{n-1}\right)  ^{2}$ for $n\geq1$. \ ML\ allows calculation
of
\[
\tilde{g}_{15}(x)=-2\lim_{n\rightarrow\infty}n^{2}\left(  \,x_{n}-\frac{1}%
{2n}+\frac{3}{8}\frac{\ln(n)}{n^{2}}\right)
\]
to high precision and \cite{F3-zeal}
\begin{align*}
\tilde{g}_{15}\left(  \frac{1}{2}\right)   &
=2.05147406138650942317849754883637246230944489958205\backslash\\
&  \;\;\;\;\;\;\;\;11023132647875850901664278100401729946190648322822...,
\end{align*}%
\begin{align*}
\tilde{g}_{15}\left(  \frac{1}{3}\right)   &
=1.57467224586735424010327033762655336904223882231552\backslash\\
&  \;\;\;\;\;\;\;\;12559730614628877089577841074745751944760863233997....
\end{align*}
We conjecture that $\delta_{15}=-(3/4)\ln(2)$ (identical to $\delta_{13}$) and
this holds for any choice of $x=x_{0}$. \ More generally,
\[
\delta(k)=-\frac{k+1}{2k}\ln(k)
\]
corresponding to $\theta(x)=x\left(  1-x\right)  ^{k}$ for integer $k\geq1$.
\ A graph of $\tilde{g}_{15}\left(  x\right)  $ appears in \cite{F4-zeal},
associated with $x(1-a\,x+x^{2})$ for $a=2$. \ The progression of graphs for
$a=3/2$, $5/3$, $\sqrt{3}$, $9/5$ and $19/10$ is quite intriguing and deserves
to be better known. \ 

We finally turn to $\theta_{16}(x)=x(1-3x+3x^{2})$, for which $\tau=1$ and
$\gamma=-3$. \ EJ\ gives%
\begin{align*}
g_{16}(x)  &  =\frac{1}{3x}+\frac{2}{3}\ln(x)+\frac{2}{3}x+\frac{5}{6}%
x^{2}+\frac{14}{9}x^{3}+\frac{89}{30}x^{4}+\frac{563}{150}x^{5}-\frac
{8749}{1260}x^{6}-\frac{22259}{294}x^{7}\\
&  -\frac{128419}{420}x^{8}-\frac{4596701}{9450}x^{9}+\frac{19562713}%
{8250}x^{10}+\frac{4523425099}{254100}x^{11}+\frac{495227461}{51480}x^{12}\\
&  -\frac{6781873063663}{13663650}x^{13}-\frac{171780306415997}{73573500}%
x^{14}+\frac{1764254963282093}{157657500}x^{15}+\cdots.
\end{align*}
To assess $g_{16}(x)$ for some $0<x<1$, define $x_{0}=x$ and $x_{n}%
=x_{n-1}\left(  1-3x_{n-1}+3x_{n-1}^{2}\right)  $ for $n\geq1$. \ ML\ allows
calculation of%
\[
\tilde{g}_{16}(x)=-3\lim_{n\rightarrow\infty}n^{2}\left(  \,x_{n}-\frac{1}%
{3n}+\frac{2}{9}\frac{\ln(n)}{n^{2}}\right)
\]
to high precision and \cite{F3-zeal}%
\begin{align*}
\tilde{g}_{16}\left(  \frac{1}{2}\right)   &
=1.11293077809135353800166857803390442442494189554367\backslash\\
&  \;\;\;\;\;\;\;\;24343525392480517435909954503132620225124156973073...,
\end{align*}%
\begin{align*}
\tilde{g}_{16}\left(  \frac{1}{3}\right)   &
=1.35456732398262548307271488105822274526124339313976\backslash\\
&  \;\;\;\;\;\;\;\;49107656609639544183939435929936110800486434074608....
\end{align*}
We conjecture that $\delta_{16}=-(2/3)\ln(3)$ and this holds for any choice of
$x=x_{0}$. \ Formal proofs would be good to see someday. \ Updates will be
reported in the Addenda.%

\begin{figure}
[ptb]
\begin{center}
\includegraphics[
height=3.3053in,
width=3.1704in
]%
{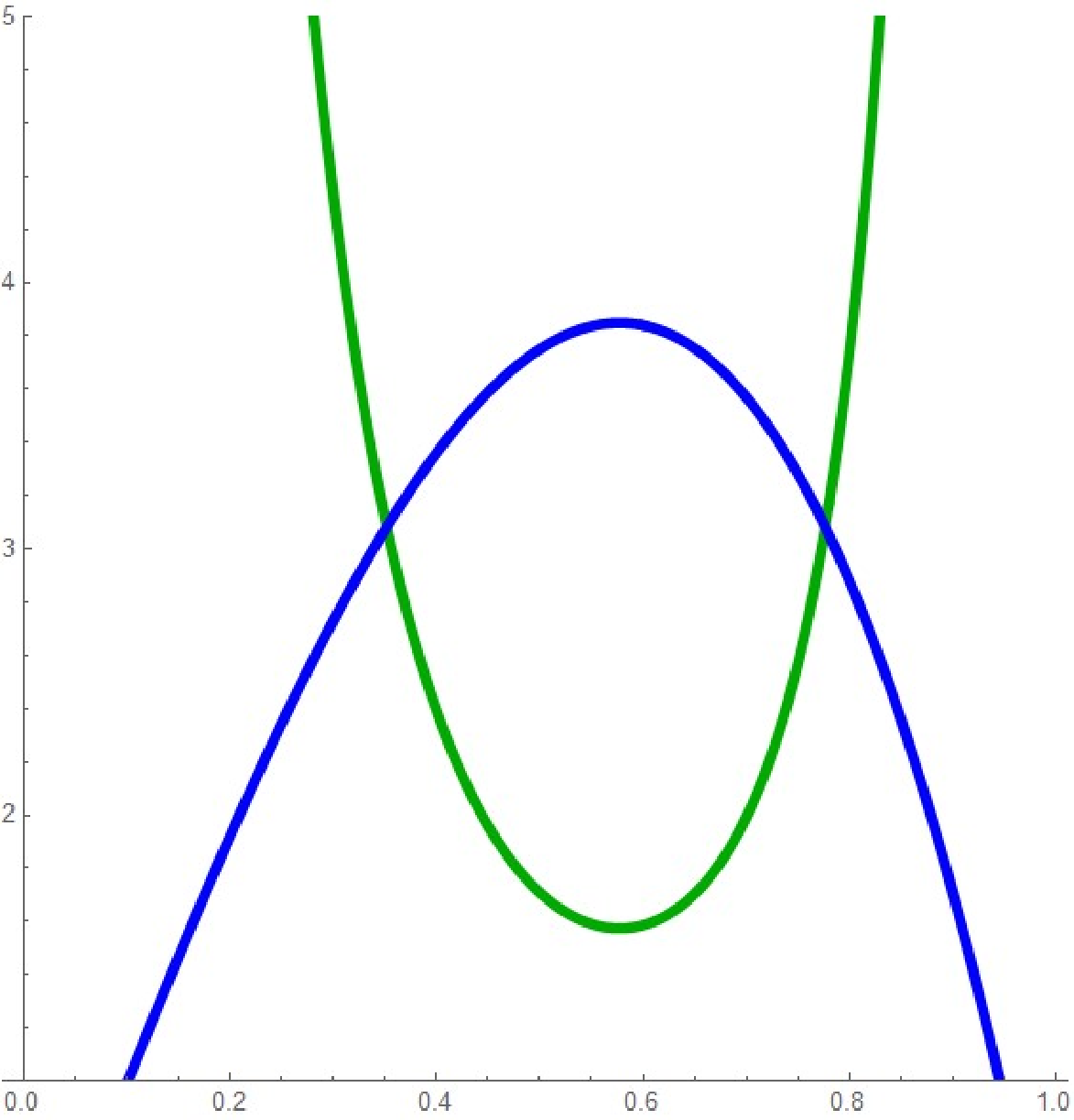}%
\caption{Blue curve is $10\theta_{13}(x)$, scaled for visibility; green curve
is ML-based $\tilde{g}_{13}(x)$. Lower left-hand corner is $(0,1)$. Distance
between vertical notches is eight times that for horizontal.}%
\end{center}
\end{figure}
\begin{figure}
[ptb]
\begin{center}
\includegraphics[
height=2.8496in,
width=5.4708in
]%
{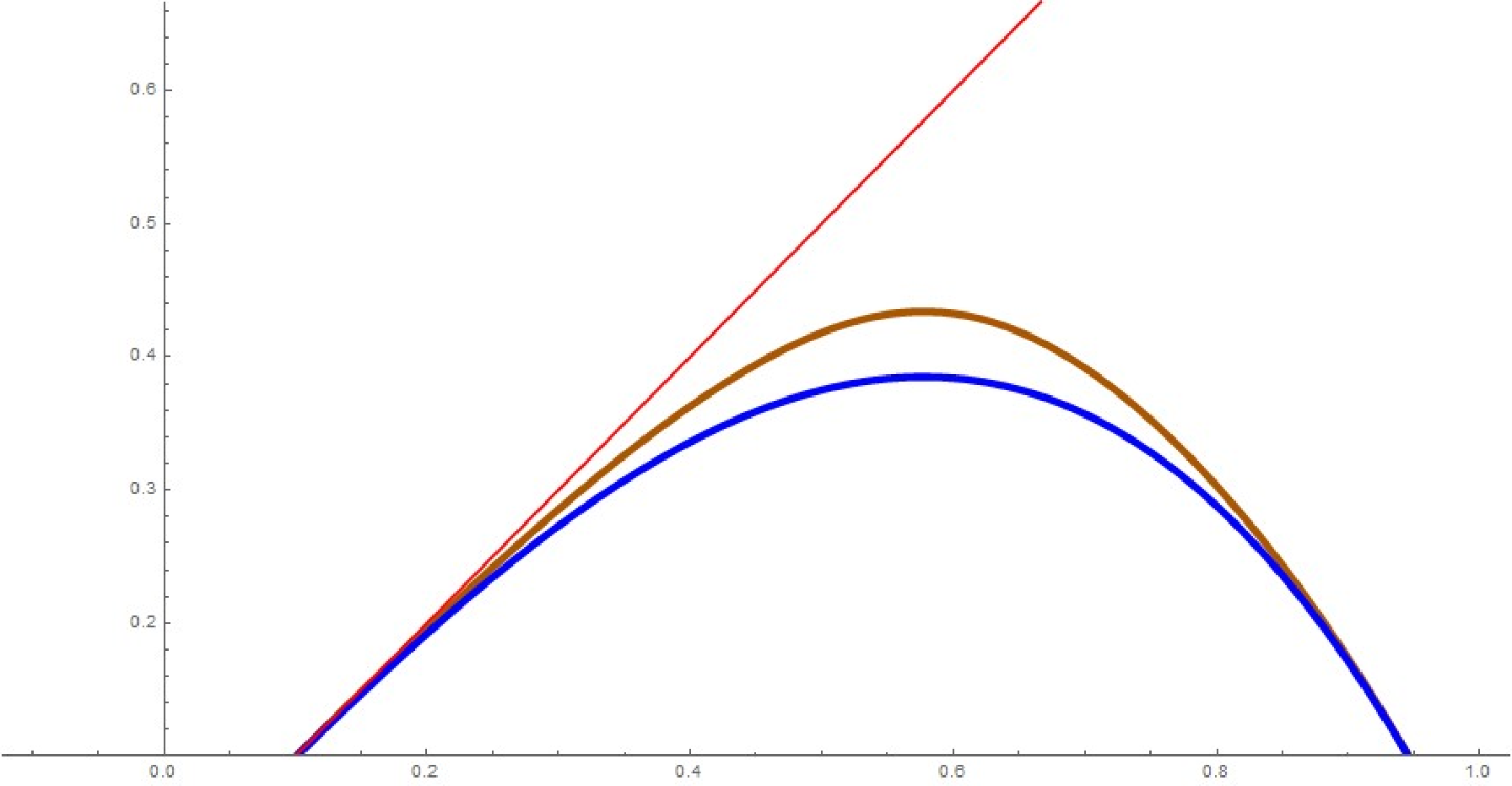}%
\caption{Half-iterate of $x\left(  1-x^{2}\right)  $ in orange; $x\left(
1-x^{2}\right)  $ in blue; diagonal line in red. \ Lower left-hand corner is
$(0,1/10)$. }%
\end{center}
\end{figure}

\section{Clarification}

A software programming error led to false conclusions in an earlier draft
(v3). \ The pseudocode in \cite{F1-zeal, F2-zeal} is entirely correct, but
(inexplicably, over time) the expression\ $%
{\textstyle\sum\nolimits_{m=2}^{K-2}}
v_{m}x^{\tau\,m+1}$ became altered to $%
{\textstyle\sum\nolimits_{m=3}^{K-1}}
v_{m}x^{\tau\,m-1}$. \ We mention this for completeness' sake; every EJ-based
result in that draft about $g_{14}$, $g_{15}$, $g_{16}$ is regrettably wrong.

\section{Addendum I}

The quartic recurrence%
\[%
\begin{array}
[c]{ccc}%
\xi_{n}=\dfrac{3}{4}+\dfrac{1}{4}\xi_{n-1}^{4}, &  & 0\leq\xi_{0}<1
\end{array}
\]
arises in the study of random\ Galton-Watson quaternary tree heights at
criticality.\ \ Letting $\xi_{n}=1-4x_{n}$, we obtain%
\begin{align*}
1-4x_{n}  &  =\frac{1}{4}\left[  3+\left(  1-4x_{n-1}\right)  ^{4}\right] \\
&  =\frac{1}{4}\left(  4-16x_{n-1}+96x_{n-1}^{2}-256x_{n-1}^{3}+256x_{n-1}%
^{4}\right) \\
&  =1-4x_{n-1}+24x_{n-1}^{2}-64x_{n-1}^{3}+64x_{n-1}^{3}%
\end{align*}
hence%
\[%
\begin{array}
[c]{ccc}%
x_{n}=\theta_{17}(x_{n-1})=x_{n-1}\left(  1-6x_{n-1}+16x_{n-1}^{2}%
-16x_{n-1}^{3}\right)  , &  & 0<x_{0}<1/2.
\end{array}
\]
ML\ allows calculation of%
\[
\tilde{g}_{17}(x)=-6\lim_{n\rightarrow\infty}n^{2}\left(  \,x_{n}-\frac{1}%
{6n}+\frac{5}{54}\frac{\ln(n)}{n^{2}}\right)
\]
to high precision and
\begin{align*}
\tilde{g}_{17}\left(  \frac{1}{4}\right)   &
=1.18801384280361604567470036016221695576085770053662\backslash\\
&  \;\;\;\;\;\;\;\;68477703163061637332612940383107769158426699380380....
\end{align*}
We conjecture that $\delta_{17}=-(5/9)\ln(6)$ and this holds for any choice of
$x=x_{0}$. \ 

\section{Addendum II}

The cubics in this paper are all of type $x\left(  1-\sigma\,x-\rho
\,x^{2}\right)  $, where $-\infty<\rho<\infty$ and $\sigma\geq0$. \ We
conjecture, assuming $\sigma>0$, that%
\[
\delta=\left(  1+\frac{\rho}{\sigma^{2}}\right)  \ln\left(  \frac{1}{\sigma
}\right)  .
\]
This extends an earlier formula connected with $0<\rho<1$ and $\sigma=1-\rho$. \ 

A\ more far-reaching conjecture starts with observations that both%
\[%
\begin{array}
[c]{ccccc}%
x(1-x)^{3}=x\left(  1-3x+3x^{2}-x^{3}\right)  &  & \text{and} &  & \theta
_{16}(x)=x\left(  1-3x+3x^{2}\right)
\end{array}
\]
share $\delta=-(2/3)\ln(3)$; both%
\[%
\begin{array}
[c]{ccccc}%
\theta_{17}(x)=x\left(  1-6x+16x^{2}-16x^{3}\right)  &  & \text{and} &  &
x\left(  1-6x+16x^{2}\right)
\end{array}
\]
share $\delta=-(5/9)\ln(6)$; both \cite{F1-zeal}%
\[%
\begin{array}
[c]{ccccc}%
\theta_{5}(x)=\ln(1+x) &  & \text{and} &  & x\left(  1-\dfrac{1}{2}x+\dfrac
{1}{3}x^{2}\right)
\end{array}
\]
share $\delta=-(1/3)\ln(2)$; both \cite{F2-zeal}%
\[%
\begin{array}
[c]{ccccc}%
\theta_{6}(x)=x\exp(-x) &  & \text{and} &  & x\left(  1-x+\dfrac{1}{2}%
x^{2}\right)
\end{array}
\]
share $\delta=0$; both \cite{F2-zeal}%
\[%
\begin{array}
[c]{ccccc}%
\theta_{8}(x)=\dfrac{x}{1+x^{2}} &  & \text{and} &  & x\left(  1-x^{2}%
+x^{4}\right)
\end{array}
\]
share $\delta=-(1/4)\ln(2)$; both \cite{F1-zeal}%
\[%
\begin{array}
[c]{ccccc}%
\theta_{3}(x)=\sin(x) &  & \text{and} &  & x\left(  1-\dfrac{1}{6}x^{2}%
+\dfrac{1}{120}x^{4}\right)
\end{array}
\]
share $\delta=(3/5)\ln(3)$. \ In essence, $\delta$ seems to be completely
determined by the parameter $\tau$ and two coefficients $c_{1}$ \&\ $c_{2}$. \ 

\section{Addendum III}

Given a polynomial of type $x\left(  1-\sigma\,x^{\tau}-\rho\,x^{2\tau
}\right)  $ where $\sigma>0$, we believe that
\[
\delta=\left(  \frac{\tau+1}{2\tau}+\frac{\rho}{\tau\,\sigma^{2}}\right)
\ln\left(  \frac{1}{\tau\,\sigma}\right)  .
\]
This extends to transcendental functions, e.g., the Fresnel cosine
\cite{F0-zeal} has $\delta=\frac{55}{108}\ln\left(  \frac{10}{\pi^{2}}\right)
$. \ 

\section{Acknowledgements}

I am grateful to Daniel Lichtblau at Wolfram Research for kindly answering my
questions, e.g., about generalizing my original Mathematica code for ML to
arbitrary $k$. \ William Jagy \cite{J1-zeal} assisted me in more ways than he
can imagine. \ Dmitrii Kouznetsov \cite{K1-zeal} fleetingly mentioned both
$g_{3}(\pi/2)$ and $\tilde{g}_{3}(\pi/2)$; see \cite{F1-zeal}. \ The creators
of Mathematica earn my gratitude every day:\ this paper could not have
otherwise been written. \ An interactive computational notebook about ML is
available \cite{F5-zeal} which might be useful to interested readers; an
analog\ for EJ is forthcoming.

\end{document}